\documentclass[11pt,leqno]{article}
\usepackage{amscd,amsfonts,amsmath,amssymb,amsthm,amstext,latexsym}
\usepackage{enumerate,pifont,epsf
  ,graphicx}
\usepackage{color}
\usepackage[english]{babel}
\usepackage[T1]{fontenc}
\usepackage{epsfig}
\usepackage[latin1]{inputenc}

\newcommand{\R}{{\mathbb R}}
\newcommand{\s}{{\mathbb S}}
\newcommand{\C}{{\mathbb C}}
\newcommand{\N}{{\mathbb N}}
\newcommand{\0}{{\bf 0}}
\def\H{\mathbb H}

\begin{document}

\begin{title} {Non-proper complete minimal surfaces embedded in
    $\H^2\times\R$}
\end{title}
\vskip .2in

\begin{author} {Magdalena Rodr\'\i guez\thanks{Research partially
      supported by the MCyT-Feder research project MTM2007-61775 and
      the Regional J. Andaluc\'\i a Grant no. P09-FQM-5088.} and
    Giuseppe Tinaglia~\thanks{Partially supported by EPSRC grant
      no. EP/I01294X/1}}
\end{author}

\date{}
\maketitle

\begin{abstract}
  Examples of complete minimal surfaces properly embedded in
  $\H^2\times \R$ have been extensively studied and the literature
  contains a plethora of nontrivial ones.  In this paper we construct
  a large class of examples of complete minimal surfaces embedded in
  $\H^2\times\R$, not necessarily proper, which are invariant by a
  vertical translation or by a hyperbolic or parabolic screw
  motion. In particular, we construct a large family of non-proper
  complete minimal disks embedded in $\H^2\times\R$ invariant by a
  vertical translation and a hyperbolic screw motion and whose
  importance is twofold. They have finite total curvature in the
  quotient of $\H^2\times\R$ by the isometry, thus highlighting a
  different behaviour from minimal surfaces embedded in $\R^3$
  satisfying the same properties. They show that the Calabi-Yau
  conjectures do not hold for embedded minimal surfaces in $\H^2\times
  \R$.
\end{abstract}

\noindent{\it Mathematics Subject Classification:} Primary 53A10,
   Secondary 49Q05, 53C42

\section{Introduction}
 
Examples of complete minimal surfaces properly embedded in $\H^2\times
\R$ have been extensively studied and the literature contains a
plethora of nontrivial ones. In this paper we focus on complete
embedded examples, not necessarily proper, which are invariant by
either a vertical translation or a hyperbolic or parabolic screw
motion. Some examples with these properties, but all of them properly
embedded, have been constructed in~\cite{nr,mo,o,st,s,moro1,marr2,me}.

The key examples contained in this paper are complete minimal disks
embedded in $\H^2\times\R$ that are non-proper and invariant by a
vertical translation and a hyperbolic screw motion, we call them
helicoidal-Scherk examples.  The importance of such helicoidal-Scherk
examples is twofold in understanding the behaviour of minimal surfaces
in $\H^2\times \R$.

In addition to being non-proper, a significant feature of these
examples is that they have finite total curvature in the quotient of
$\H^2\times \R$ by the vertical translation or the hyperbolic screw
motion.  In~\cite{to} Toubiana proved that a complete embedded minimal
annulus with finite total curvature in the quotient of $\R^3$ by a
translation must the quotient of a helicoid. In~\cite{mr1} Meeks and
Rosenberg proved that Toubiana's result holds if the translation is
replaced by a screw-motion. Moreover, in the same paper they also show
that a complete embedded minimal surface with finite total curvature
in the quotient of $\R^3$ by a translation or a screw-motion must be
proper.  Our examples highlight a much different behaviour in
$\H^2\times\R$. Recently, Collin, Hauswirth and Rosenberg have studied
the conformal type and the geometry of the ends of properly embedded
minimal surfaces with finite total curvature in the quotient of
$\H^2\times\R$ by a vertical translation~\cite{chr}.  Our main
examples are related to but not included in their study.

The same examples are also of interest in relation to the Calabi-Yau
conjectures for embedded minimal surfaces~\cite{c,ch,y}. In~\cite{cm},
Colding and Minicozzi showed that a complete minimal surface embedded
in $\R^3$ with finite topology is proper. See~\cite{mr2} for a
generalization of their result. Our helicoidal-Scherk examples show
that Colding and Minicozzi's result does not hold in $\H^2\times
\R$. Note that in~\cite{cos} Coskunuzer has already constructed a
complete embedded disk in $\H^3$ which is not proper, thus showing
that Colding and Minicozzi's result does not generalize to $\H^3$. The
techniques that we use to construct our examples are completely
different from his.

The helicoidal-Scherk examples are constructed in the next section. In
the other sections we further generalize the construction and also
give examples of properly embedded minimal surfaces that are invariant
by a parabolic screw motion. These latter examples are included in the
study in~\cite{chr}.

\vspace{.2cm}

\noindent We would like to thank Laurent Hauswirth and Harold
Rosenberg for very helpful conversations.

\section{Helicoidal-Scherk examples}
\label{sec:helScherk}

In order to construct our examples we consider the Poincar\'e disk
model of $\H^2$; i.e.
\[
\H^2=\{z\in\C\ |\ |z|<1\},
\]
with the hyperbolic metric 
\[
g_{-1}=\frac{4}{(1-|z|^2)^2}|dz|^2.
\]
We denote by $\partial_\infty\H$ the boundary at infinity of $\H^2$
and by $\0$ the origin of $\H^2$. We use $t$ for the coordinate in
$\R$.  Finally, given any two points
$p,q\in\H^2\cup\partial_\infty\H^2$, we will denote by $\overline{pq}$
the geodesic arc joining them.

\begin{figure}
  \begin{center}
    \epsfysize=6cm \epsffile{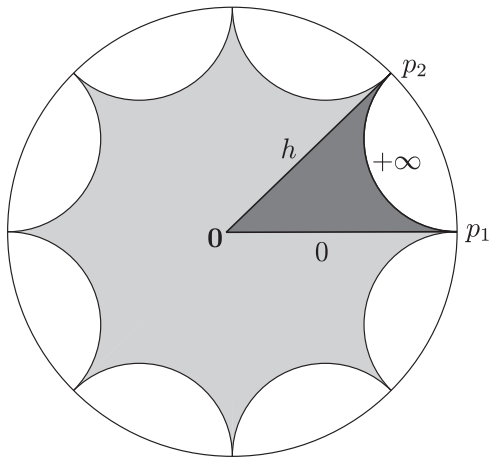}
  \end{center}
  \caption{Fundamental piece of a helicoidal-Scherk example for n=2.}
  \label{Fig1}
\end{figure}

Let us consider $p_{1}=1$ and $p_{2}=e^{i\frac{\pi}{2n}}$, for some
$n\in\N$. Let $\Omega$ be the region bounded by the ideal geodesic
triangle with vertices $\0, p_1$ and $p_2$ and edges $\overline{\0
  p_{1}}$, $\overline{\0 p_{2}}$ and $\overline{p_1 p_2}$. By
Theorem~4.9 in~\cite{marr1}, there exists a minimal graph over
$\Omega$ with boundary values $0$ over $\overline{\0 p_{1}}$, $h$ over
$\overline{\0 p_{2}}$ and $+\infty$ over $\overline{p_1 p_2}$, for any
constant $h>0$ (see Figure~\ref{Fig1}).  We call this graph the
fundamental piece. Using Schwarz reflection principle, after
considering successive symmetries with respect to the horizontal
geodesics contained in the boundary of such a graph, namely
$\overline{\0 p_{1}}\times\{0\}, \overline{\0
  p_{2}}\times\{h\}\subset\H^2\times\R$, we obtain a simply-connected
minimal surface $\widehat{M}_{n,h}$ with boundary the vertical line
$\{\0\}\times\R$ and invariant by the vertical translation $T$ by
$(\0,4nh)$ and by the hyperbolic screw motion $S$ obtained by
composing the rotation by angle $\frac\pi n$ around $\0$ and the
vertical translation by $(\0,2h)$. After reflecting across the line
$\{\0\}\times \R$, we obtain a simply-connected complete minimal
surface $M_{n,h}$ invariant by both the vertical translation $T$ and
the hyperbolic screw motion $S$. We call these surfaces {\em
  helicoidal-Scherk examples}.

We observe that if we consider $h=0$ in this construction, we obtain a
Scherk graph over a symmetric ideal polygonal domain,
see~\cite{cor2,nr}.

As a  consequence of the Gauss-Bonnet Theorem applied
on a fundamental piece (see~\cite[page 1896]{cor2} for a similar
argument), $M_{n,h}$ has finite total curvature in its
quotient by both $T$ and $S$.

We next show that $M_{n,h}$ is embedded.
Let us denote by $\Omega_i$, $i=1,\cdots,4n$, the domain obtained by
rotating $\Omega$ around the origin by an angle $\frac{\pi}{2n}(i-1)$
so that $\Omega_1=\Omega$, and let
$\widetilde\Omega_i=\overline\Omega_i\times\R$, where
$\overline\Omega_i$ is the closure of $\Omega_i$. Note that the domain
$\Omega_{2n+i}$ can also be obtained by reflecting $\Omega_i$ across
the origin. We are going to prove that
\[
M_{n,h}\cap \widetilde\Omega_{2n +1}
\] 
has no self-intersections. After this, repeating the same argument
shows that $M_{n,h}$ is embedded.  Let $p^*_{1},p^*_{2}$ be the
reflection across the origin of $p_{1}$ and $ p_{2}$.  Recall that
$\widehat{M}_{n,h}\cap\widetilde \Omega_1$ consists of a graph with
boundary values $0$ over $\overline{\0 p_{1}}$, $h$ over $\overline{\0
  p_{2}}$ and $+\infty$ over $\overline{p_1 p_2}$, together with its
vertical translates by the vector $k(\0,4nh)$, $k\in\mathbb Z$.  By
construction, since we have reflected an even number of times,
$\widehat{M}_{n,h}\cap\widetilde\Omega_{2n+1}$ consists of a union of
graphs with boundary values
\[
\left\{\begin{array}{ll}
    2nh+4knh     & \mbox{, over } \overline{\0 p^*_{1}}\\
    (2n+1)h+4knh & \mbox{, over }  \overline{\0 p^*_{2}}\\
    +\infty & \mbox{, over } \overline{p^*_1 p^*_2}
  \end{array}\right.
\]
with $k\in \mathbb Z$.  The reflection of
$\widehat{M}_{n,h}\cap\widetilde\Omega_1$ across $\{\0\}\times\R$
instead consists of a union of graphs with boundary values
\[
\left\{\begin{array}{ll}
    4knh     & \mbox{, over } \overline{\0 p^*_{1}}\\
    h+4knh & \mbox{, over }  \overline{\0 p^*_{2}}\\
    +\infty & \mbox{, over } \overline{p^*_1 p^*_2}
  \end{array}\right.
\]
with $k\in \mathbb Z$.  In fact, $M_{n,h}\cap \widetilde\Omega_{2n
  +1}$ consists of the reflected fundamental piece together with its
vertical translates by the vector $k(\0,2nh)$, $k\in\mathbb Z$. In
particular $M_{n,h}\cap \widetilde\Omega_{2n +1}$ is embedded.
Repeating this argument proves that $M_{n,h}$ is embedded.

Observe that the previous argument also shows that we can consider the
quotient of $M_{n,h}$ by $(\0,2nh)$, obtaining a non-orientable
complete non-proper embedded minimal surface.

Finally we remark that $M_{n,h}\cap\widetilde\Omega_1$ accumulates to
$\overline{p_1 p_2}\times\R$, and therefore $M_{n,h}$ is a
simply-connected (in particular, with finite topology) minimal surface
embedded in $\H^2\times\R$ which is complete but not proper.
\vspace{.5cm}

Let us now describe a generalization of these examples. Instead of
considering a geodesic triangle, let $\Omega$ be the region bounded by
an ideal geodesic polygon constructed in the following way. As before,
let $p_{1}=1$ and $p_{2}=e^{i\frac{\pi}{2n}}$ and let
$\text{arc}(p_1p_2)$ denote the shortest arc in $\partial_\infty\H^2$
with end points $p_1,p_2$.  In this construction, the geodesics
$\overline{\0 p_{1}}$ and $\overline{\0 p_{2}}$ are the same but,
instead of connecting the two with the geodesic $\overline{p_1 p_2}$,
we consider $k\geq 1$ points $q_1,\cdots,q_k$, cyclically ordered in
$\text{arc}(p_1p_2)$. We define $\Omega$ as the region bounded by the
ideal geodesic polygon with vertices $\0, p_1, q_1,\cdots,q_k$ and
$p_2$. Assuming that $\Omega$ satisfies the Jenkins-Serrin condition
of Theorem~4.9 in~\cite{marr1}, we can find a graph over $\Omega$ with
boundary values $0$ on $\overline{\0 p_1}$, $h>0$ on $\overline{\0
  p_2}$ and alternating $\pm\infty$ on the remaining geodesic arcs
$\overline{p_1q_1}$, $\overline{q_1 q_2}$, \ldots, $\overline{q_{k-1}
  q_k}$, $\overline{q_k p_2}$. After considering successive symmetries
with respect to the horizontal and vertical geodesics contained in the
boundary of such a graph, we obtain a simply-connected complete
embedded minimal surface in $\H^2\times \R$ invariant by the same
vertical translation $T$ and the same hyperbolic screw motion $S$
previously defined. Again, this surface is non-proper and has finite
total curvature when considered in the quotient by $T$ or
$S$. Moreover, it admits a non-orientable quotient by the vertical
translation given by the vector $(\0,2nh)$.  We also refer to such
examples as {\em helicoidal-Scherk examples}.

It is easy to show that the class of these more general examples is
rather large. Here is an easy way to construct domains as previously
described.  If for $j=1,\dots, k$, we let
$q_j=e^{i\frac{j\pi}{2n(k-1)}}$, then we recover the symmetric
helicoidal-Scherk examples for a smaller choice of $h$, by the
generalized maximum principle for such minimal graphs
(see~\cite[Theorem 2]{cor2} or~\cite[Theorem 4.13]{marr1}).  However,
after slightly perturbing one such $q_i$, we would obtain a domain
satisfying the Jenkins-Serrin condition of Theorem~4.9
in~\cite{marr1}.  In particular, we observe that when $k=1$ and $q_1$
is any point in $\text{arc}(p_1p_2)$ then the Jenkins-Serrin condition
is satisfied.\vspace{.5cm}

As mentioned in the introduction, the importance of these examples is
twofold.

\begin{itemize}
\item In~\cite{to} Toubiana proved that a complete embedded minimal
  annulus with finite total curvature in the quotient of $\R^3$ by a
  translation must the quotient of a helicoid. In~\cite{mr1} Meeks and
  Rosenberg proved that Toubiana's result holds if the translation is
  replaced by a screw-motion. Moreover, in the same paper they also
  show that a complete embedded minimal surface with finite total
  curvature in the quotient of $\R^3$ by a translation or a
  screw-motion must be proper.  Our examples highlight a much
  different behaviour in $\H^2\times\R$.

\item In~\cite{cm}, Colding and Minicozzi showed that a complete
  minimal surface embedded in $\R^3$ with finite topology is
  proper. Thus showing that the Calabi-Yau conjectures hold for
  complete minimal surfaces embedded in $\R^3$, see~\cite{c,ch,y}. Our
  helicoidal-Scherk examples show that Colding and Minicozzi's result
  does not hold in $\H^2\times \R$. Note that in~\cite{cos} Coskunuzer
  has already constructed a complete embedded disk in $\H^3$ which is
  not proper, thus showing that Colding and Minicozzi's result does
  not generalize to $\H^3$. The techniques that we have used to
  construct the helicoidal-Scherk examples are completely different
  from his. Note also that in~\cite{mr2}, Meeks and Rosenberg
  generalized the result in~\cite{cm} to complete minimal surfaces
  with positive injectivity radius and, among other things, showed
  that the closure of a complete minimal surface with positive
  injectivity radius embedded in a 3-manifold has the structure of a
  minimal lamination. The closure of a helicoidal-Scherk example is
  the minimal lamination given by the union of such helicoidal-Scherk
  example with the related totally geodesic vertical
  planes. \end{itemize}

\section{Helicoidal examples}
\label{sec:hiphel}

\begin{figure}
  \begin{center}
    \epsfysize=6cm \epsffile{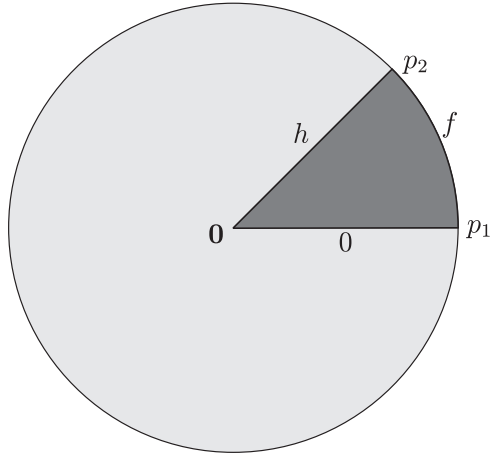}
  \end{center}
  \caption{Fundamental
    piece of a helicoidal example with $m=4$.}
  \label{Fig2}
\end{figure}

Let us now consider $p_{1}=1$ and $p_{2}=e^{i\frac{\pi}{m}}$, with
$m\in\N$.  Let $\Omega$ be the region bounded by $\overline{\0 p_1},
\overline{\0 p_2}$ and $\text{arc}(p_1p_2)$, see Figure~\ref{Fig2}. By
Theorem 4.9 in \cite{marr1}, there exists a minimal graph over
$\Omega$ with boundary values $0$ over $\overline{\0 p_1}$, $h$ over
$\overline{\0 p_2}$ and $f$ over $\text{arc}(p_1p_2)$, for any $h>0$
and any continuous function $f$ on $\text{arc}(p_1p_2)$ (in fact,
finitely many points of discontinuity for $f$ are allowed). Again,
after considering successive symmetries with respect to the horizontal
geodesics contained in the boundary of such a graph, we get a minimal
surface $\widehat{M}_{m,h,f}$ bounded by the vertical line
$\{\0\}\times\R$ and invariant by the vertical translation $T$ by
$(\0,2mh)$ and by the hyperbolic screw motion $S$ obtained by
composition of the rotation by angle $\frac{2\pi}m$ around $\0$ and
the vertical translation by $(\0,2h)$. Considering a final symmetry
with respect to $\{\0\}\times\R$, we obtain a simply-connected
complete minimal surface $M_{m,h,f}\subset\H^2\times \R$ which is
invariant by the vertical translation $T$ and by the hyperbolic screw
motion $S$. We call these surfaces {\em helicoidal examples}. These
examples have infinite total curvature in the quotient, since their
normal vectors do not become horizontal when we approach points in
$\text{arc}(p_1p_2)\times\R$ (see~\cite[Theorem 3.1]{hr}).

Let us show that $M_{m,h,f}$ is embedded when $m$ is even and, if $f$
satisfies certain conditions, when $m$ is odd.  Using the same
notation as in the previous section, we know that
$\widehat{M}_{m,h,f}\cap\widetilde\Omega_{m+1}$ consists of the union
of the minimal graph with boundary values
\[
\left\{\begin{array}{ll}
    mh     & \mbox{, over } \overline{\0 p^*_{1}}\\
    (m+1)h & \mbox{, over }  \overline{\0 p^*_{2}}\\
    f_m & \mbox{, over } \text{arc}(p_1^*p_2^*)
  \end{array}\right.
\]
where
\[
f_m= \left\{\begin{array}{ll}
    mh+f     & \mbox{, if } m  \mbox{ is even}\\
    (1+m)h-f & \mbox{, if } m \mbox{ is odd}
  \end{array}\right.
\]
together with its vertical translates by the vector $k(\0,2mh)$,
$k\in\mathbb Z$.  The reflection of
$\widehat{M}_{m,h,f}\cap\widetilde\Omega_1$ across $\{\0\}\times\R$
instead consists of the graph with boundary values
\[
\left\{\begin{array}{ll}
    0     & \mbox{, over } \overline{\0 p^*_{1}}\\
    h & \mbox{, over }  \overline{\0 p^*_{2}}\\
    f & \mbox{, over } \text{arc}(p_1^*p_2^*)
  \end{array}\right.
\]
together with its vertical translates by the vector $k(\0,2mh)$,
$k\in\mathbb Z$.  Hence, using the general maximum principle for
minimal graphs~\cite[Theorem 4.16]{marr1}, we get that $M_{m,h,f}\cap
\widetilde\Omega_{m +1}$ is embedded when $m$ is even or when $m$ is
odd and
\[
(1-m)h\leq 2f\leq (1+m)h .
\]
By symmetry, $M_{m,h,f}$ is embedded under the same conditions.

From the argument above we deduce that, when $m$ is even,
$M_{m,h,f}\cap \widetilde\Omega_{m+1}$ consists of the reflected
fundamental piece together with its vertical translates by the vector
$k(\0,mh)$, $k\in\mathbb Z$. Thus $M_{m,h,f}$ admits also a
non-orientable quotient by $(\0,mh)$.

In the case $m$ is even, if we consider the sequence of functions
$\{f_k\}$, where $f_k=k$ over $\text{arc}(p_1p_2)$, then we obtain the
fundamental piece of the corresponding symmetric helicoidal-Scherk
example as a limit of the fundamental piece of $M_{m,h,f_k}$ and thus
$M_{\frac m2,h}$ as a limit of the sequence of surfaces
$\{M_{m,h,f_k}\}_k$. We could take another choice of functions $f_k$
with the same limit but in such a way that each $M_{m,h,f_k}$ has a
smooth boundary.  In fact, any helicoidal-Scherk example can be
recovered as a limit of some sequence $\{M_{m,h,f_k}\}_k$ of
helicoidal examples, by choosing appropriate functions $f_k$.

Finally, observe that if we consider $f(e^{i\, t})= \frac{hm}{\pi}\,
t$, with $t\in(0,\frac\pi m)$, we recover one of the helicoids given
by Nelli and Rosenberg in~\cite{nr}, congruent to the Euclidean
one. In fact, by varying $h>0$ we re-obtain all of their examples.
Hence the family of helicoidal examples contains the helicoids.

\section{Helicoidal-Scherk examples with axis at infinity}

\begin{figure}
  \begin{center}
    \epsfysize=6cm \epsffile{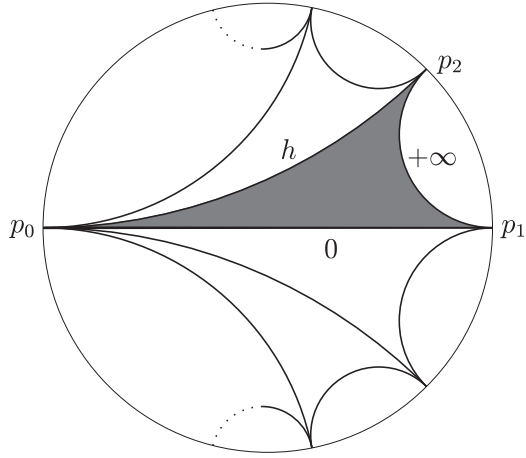}
  \end{center}
  \caption{Fundamental piece of a helicoidal-Scherk example with axis
    at infinity.}
  \label{Fig3}
\end{figure}

We now take $p_0=-1$, $p_{1}=1$ and $p_{2}=e^{i\theta}$, for some
$\theta\in(0,\pi)$. Let $\Omega$ be the region bounded by the ideal
geodesic triangle with vertices $p_0, p_1$ and $p_2$, see
Figure~\ref{Fig3}. By Theorem 4.9 in~\cite{marr1}, there exists a
minimal graph over $\Omega$ with boundary values $0$ over
$\overline{p_0 p_{1}}$, $h$ over $\overline{p_0 p_{2}}$ and $+\infty$
over $\overline{p_1 p_2}$, for any constant $h>0$.  After considering
successive symmetries with respect to $\overline{p_0
  p_{1}}\times\{0\}, \overline{p_0
  p_{2}}\times\{h\}\subset\H^2\times\R$, we obtain a properly embedded
minimal surface ${M}_{\theta,h}$ (in fact, it is a graph over an ideal
polygonal domain with infinitely many boundary geodesic arcs)
invariant by the parabolic screw motion $P$ obtained by composition of
the parabolic translation with fixed point $p_0$ which maps $p_1$ onto
$e^{i2\theta}$ with the vertical translation by $(\0,2h)$. We call
these examples {\em helicoidal.Scherk examples with axis at infinity},
since they can be obtained as a limit of helicoidal-Scherk examples
whose axes go to infinity. As a consequence of the Gauss-Bonnet
Theorem, $M_{\theta,h}$ has finite total curvature in its quotient by
$P$.

We observe that if we consider $h=0$ in this construction, we obtain a
pseudo-Scherk graph considered by Leguil and Rosenberg in~\cite{lr}.

Just like in section~\ref{sec:helScherk}, these examples can be
generalized by taking an ideal geodesic polygon $\Omega$ with vertices
$p_0=-1$, $p_{1}=1$, $p_{2}=e^{i\theta}$ and $k\geq 1$ points
$q_1,\cdots,q_k$ in $\text{arc}(p_1p_2)$, such that $\Omega$ satisfies
the Jenkins-Serrin condition of Theorem~4.9 in~\cite{marr1}. One such
polygonal domain is called pseudo-Scherk polygon in~\cite{lr}. We
start with the graph over $\Omega$ with boundary values $0$ on
$\overline{p_0 p_1}$, $h>0$ on $\overline{p_0 p_2}$ and alternating
$\pm\infty$ on the remaining geodesics. After considering successive
symmetries with respect to the horizontal geodesics contained in the
boundary of such a graph, we obtain a properly embedded minimal
surface in $\H^2\times\R$ invariant by the parabolic screw motion $P$
described above (again, it is a graph over an ideal polygonal domain
with infinitely many boundary geodesic arcs). In the quotient by $P$,
such a surface has finite total curvature.  We also refer to these
generalized surfaces as {\em helicoidal-Scherk examples with axis at
  infinity}.

\section{Helicoidal examples with axis at infinity}

\begin{figure}
  \begin{center}
    \epsfysize=6cm \epsffile{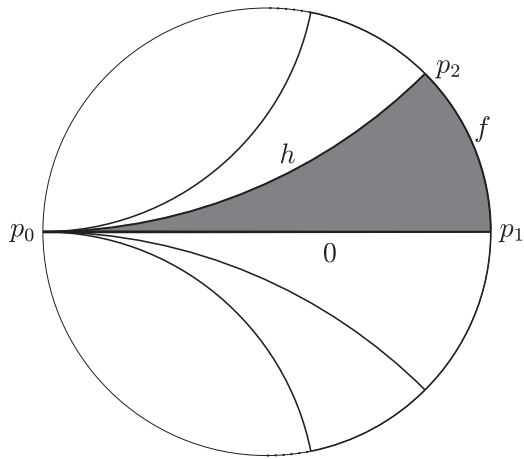}
  \end{center}
  \caption{Fundamental piece of a helicoidal example with axis at
    infinity.}
  \label{Fig4}
\end{figure}

Let us now consider $p_0=-1$, $p_{1}=1$ and $p_{2}=e^{i\theta}$, with
$\theta\in(0,\pi)$.  Let $\Omega$ be the region bounded by
$\overline{p_0 p_1}, \overline{p_0 p_2}$ and $\text{arc}(p_1p_2)$, see
Figure~\ref{Fig4}. By Theorem 4.9 in \cite{marr1}, there exists a
minimal graph over $\Omega$ with boundary values $0$ over
$\overline{p_0 p_1}$, $h$ over $\overline{p_0 p_2}$ and $f$ over
$\text{arc}(p_1p_2)$, for any $h>0$ and any continuous function $f$ on
$\text{arc}(p_1p_2)$ (again, finitely many points of discontinuity for
$f$ are allowed). After considering successive symmetries with respect
to the horizontal geodesics contained in the boundary of such a graph,
we get a properly embedded minimal surface ${M}_{\theta,h,f}$ (in
fact, it is an entire graph) invariant by the parabolic screw motion
$P$ obtained by composition of the parabolic translation with fixed
point $p_0$ which maps $p_1$ onto $e^{i2\theta}$ with the vertical
translation by $(\0,2h)$. In the quotient by $P$, ${M}_{\theta,h,f}$
has infinite total curvature. We call these surfaces {\em helicoidal
  examples with axis at infinity}.

We observe that, arguing as in section~\ref{sec:hiphel}, any
helicoidal-Scherk example with axis at infinity can be recovered as a
limit of helicoidal examples $M_{\theta,h,f_k}$, by choosing
appropriate functions $f_k$.

Finally, if we consider $f(e^{i\, t})= \frac{h}{\theta}\, t$, for any
$t\in(0,\theta)$, we recover one of the examples invariant by the
1-parametric isometry group generated by $P$, founded by
Onnis~\cite{o} 
and Sa~Earp~\cite{s} independently. 

\section{Non-periodic examples}
\label{}

In this last section, we point out how this method can be used to
construct a lot of simply-connected examples which cannot be written
as graphs. We now let $p_{1}=1$ and $p_{2}=e^{i\theta}$, for some
fixed $\theta\in(0,\pi)$, and define $\Omega$ as the domain bounded by
$\overline{\0 p_1}$, $\overline{\0 p_2}$ and $\text{arc}(p_1p_2)$.  By
Theorem 4.9 in \cite{marr1}, we know there exists a minimal graph over
$\Omega$ with boundary values $+\infty$ on $\overline{\0 p_1}$, $0$ on
$\overline{\0 p_2}$ and $f$ on $\text{arc}(p_1p_2)$, for any
continuous function $f$ (again, finitely many discontinuity points are
allowed).  By rotating such a graph by an angle $\pi$ about the
horizontal geodesic $\overline{\0 p_2}\times\{0\}$ contained in its
boundary, we obtain a minimal graph whose boundary consists of the
vertical line $\{\0\}\times\R$. After extending such a graph by
symmetry about its boundary, we obtain a properly immersed
simply-connected minimal surface. When $\theta\leq \pi/2$ or $f$ is
positive, the obtained surface is embedded. And its asymptotic
boundary curve is smooth if $f=0$ on $p_2$.

When $\theta\leq \pi/2$ and $f$ diverges to $+\infty$ at any point (or
to $\pm\infty$ alternately over a finite number of arcs contained in
$\text{arc}(p_1p_2)$, with some aditional restrictions to where the
endpoints of such arcs are placed in order to satisfy the
Jenkins-Serrin condition in the limit) we get the simply-connected
minimal examples with finite total curvature constructed by Pyo and
the first author in~\cite{pr}, called {\em twisted Scherk examples}.

\bibliographystyle{plain}

\mbox{}\\

\noindent
Magdalena Rodr\'\i guez\\
Departamento de Geometr\'\i a y Topolog\'\i a\\
Universidad de Granada\\
Fuentenueva, 18071, Granada, Spain\\
e-mail: \texttt{magdarp@ugr.es}\\

\noindent
Giuseppe Tinaglia\\
Mathematics Department\\
King's College London\\
The Strand, London WC2R 2LS, United Kingdom\\
e-mail: \texttt{giuseppe.tinaglia@kcl.ac.uk}

\end{document}